\documentclass[12pt]{article}
\usepackage{amssymb}
\usepackage{graphicx}
\usepackage{epsfig}

\def \QQ{Q \! \! \! \!\prime\;\,}

\begin{document}

\title{\bf Orbits of $\QQ^*(\sqrt{k^2m})$ under the action of the modular group $PSL(2,\mathbb {Z})$ }

\author{M. Aslam Malik \thanks{aslam.math@pu.edu.pk} and  Muhammad Riaz \thanks{mriaz.math@pu.edu.pk} \\
Faculty of Science, Department of Mathematics\\
 University of the Punjab, Lahore, Pakistan.}

\date{}
\maketitle
\begin{abstract}

We look at some ways in which coset diagrams have been used to find the orbits, number of subgroups and
structure of the finitely generated groups. In this paper we use coset diagrams and modular arithmetic to
 determine the exact number of $G$-orbits of $\QQ^*(\sqrt{p^k})$, $\QQ^*(\sqrt{2p^k})$, $\QQ^*(\sqrt{2^2p^k})$,
  and in general $\QQ^*(\sqrt{2^lp^k})$, for each  $l\geq3$ and $k=2h+1\geq3$, for each odd prime $p$.\\
\end{abstract}
AMS Mathematics subject classification (2000): 05C25, 11E04, 20G15\\
{\bf Keywords:} Finitely generated groups; Coset diagrams; Orbits.\\
\section{Introduction}
In algebra and geometry, a group action is a way of describing symmetries of objects using groups.
The essential elements of the object are described by a set and the symmetries of the object are
described by the symmetric group of this set, which consists of bijective transformations on the set.
The symmetries play an important role in the classical and quantum mechanics. It is well known that the modular group
 $PSL(2,\mathbb {Z})$ has finite presentation  $G=\langle x,y:x^2=y^3=1\rangle$  where
   $x:\alpha\rightarrow\frac{-1}{\alpha}, y:\alpha\rightarrow\frac{\alpha-1}{\alpha}$ are elliptic
transformations and their fixed points in the upper half plane are
$i$ and $e^{2\pi i/3}$.\\
Coset diagrams have been used to study quotients, orbits, number of
subgroups and structure of the finitely generated groups. The
concept of coset diagrams seems to originate from the work by
Schreier and Reidemeister in the 1920s. The interest has been risen
in the last decades as the possibilities to use the coset diagram
techniques in combination with mathematical software has been
improved. Novotny and Hrivnak \cite{Novotny}  consider the action of
the finitely generated group $SL(m,\mathbb {Z}_n)$ on the ring $
\mathbb{Z}^m_n$ and determined the orbits for $n$ arbitrary natural
number.   Shabbir and Khan \cite{Gshabbir}
 discussed conformally flat- but non flat Bianchi type I and cylindrically symmetric static space-times according
 to proper projective symmetry by using some algebraic and direct integration techniques. It is shown that the special
  class of the above space-times admit proper projective vector fields.\\
Torstensson used coset diagrams to study the quotients of the
modular group in \cite{Anna}.  The Number of Subgroups of
$PSL(2,\mathbb {Z})$ when acting on  $F_p\cup\{\infty\} $ has been
discussed in \cite{R14} and the subgroups of the classical modular
group  has been discussed in \cite{R15}.   Higman and  Mushtaq
(1983) introduced the concept of the coset diagrams for the modular
group $PSL(2,\mathbb {Z})$ and laid its foundation. Mushtaq (1988)
showed that for a fixed non-square positive integer $n$, there are
only a finite number of ambiguous numbers in $\QQ^*(\sqrt{n})$  and
that the ambiguous numbers in the coset diagram for the orbit
$\alpha^G$ form a closed path and it is the only closed path
contained in it.  By using the coset diagrams for the orbit of the
modular group $G=\langle x,y:x^2=y^3=1\rangle$ acting on  the real
quadratic fields M. A. Malik et al determined the exact number of
ambiguous numbers in $\QQ^*(\sqrt{n})$ in \cite{R11},  as a function
of $n$. The ambiguous length of an orbit $\alpha^G$ is the  number
of ambiguous numbers in the same orbit. M. A. Malik et al. in
\cite{R12} proved that $\QQ^*(\sqrt{p})$,  $p\equiv1~(mod~4)$,
splits into at least two orbits namely
 $(\sqrt{p})^G$ and $(\frac{1+\sqrt{p}}{2})^G$, and it was   also proved by the same authors that $\QQ^*(\sqrt{p})$,
 $p\equiv3~(mod~4)$, splits into at least two orbits namely  $(\sqrt{p})^G$ and $(\frac{\sqrt{p}}{-1})^G$ \cite{R5}.
In \cite{Asim} it was proved that there exist two proper $G$-subsets of $\QQ^*(\sqrt{n})$ when $n\equiv0~(mod~p)$ and four $G$-subsets
of $\QQ^*(\sqrt{n})$ when $n\equiv0~(mod~pq)$. In \cite{Riaz1} we generalized these result for $n\equiv 0~(mod~p_1 p_2...p_r)$,
where $p_1, p_2,...p_r$ are distinct odd primes. We also proved for $h=2k+1\geq3$
there are exactly two $G$-orbits of $\QQ^*(\sqrt{2^h})$ namely $ (2^k\sqrt{2})^G$ and $(\frac{2^k\sqrt{2}}{-1})^G$.\\
By a circuit $(n_k,...n_2,n_1)$   we always mean a closed path in which $n_1$ triangles have one vertex outside
   the circuit and $n_2$ triangles have one vertex inside the    circuit and so on. This circuit induces an element
   $(yx)^{n_k}...(y^2x)^{n_2}(yx)^{n_1}$ of $G$  which fixes a particular vertex $k$,
    where $k$ must be an ambiguous number, the detail is given in \cite{R10}.\\
Example: By a circuit $(2,1,3,1,2,1)$ we mean the    transformation\\
    $g=(yx)^{2}(y^2x)^{1}(yx)^{3}(y^2x)^{1}(yx)^{2}(y^2x)^{1}$  which fixes a particular vertex $k$,
    that is $(g)k=k$ as shown in figure 1. \\
\begin{center}
\includegraphics[width=4.0in, height=3.0in]{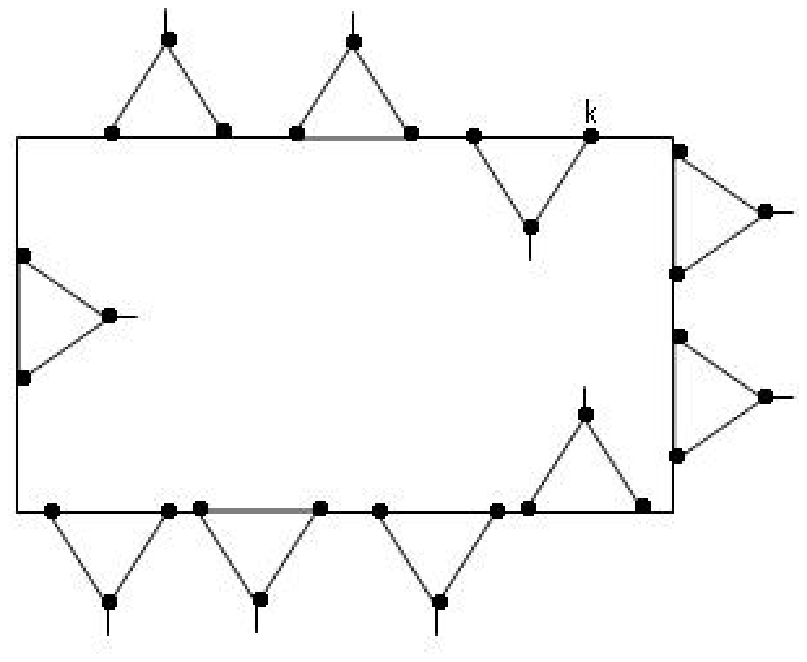}\\
\end{center}
\centerline{Fig: 1}
\textbf{Theorem 1.1}\cite{Asim}\\
Let $n\equiv 0(mod~2^3)$ and $\alpha=\frac{a+\sqrt{n}}{c}\in
\QQ^*(\sqrt{n})$ with $b=\frac{a^2-n}{c}$. Then $A=\{\alpha\in
\QQ^*(\sqrt{n}):~\textmd{either}~c\equiv i(mod~2^3)~\textmd{or
if}~c\equiv0~\textmd{or}~4(mod~2^2)~\textmd{then}~b\equiv
i(mod~2^3)\}$ is a $G$-subset of $\QQ^*(\sqrt{n})$ for each $i=1,3,5$ and $7$.\\
\\
\textbf{Theorem 1.2}\cite{R4}\\
Let  $p$ be an odd prime and $n\equiv 0(mod~p^k)$; $k\geq 1$,
Take $\alpha=\frac{a+\sqrt{n}}{c}\in \QQ^*(\sqrt{n})$ with
$b=\frac{a^2-n}{c}$. Then\\
 $A_1=\{\alpha\in \QQ^*(\sqrt{n}):(c/p)=1 ~~\textmd{or}~~(b/p)=1 \}$
 and\\
 $~A_2=\{\alpha\in \QQ^*(\sqrt{n}):(c/p)=-1 ~~\textmd{or}~~(b/p)=-1\}$
 are exactly two $G$-subsets of $Q^*(\sqrt{n})$, depending upon the classes $[a,b,c](mod~p^k)$\\

\section{ $G$-orbits of $\QQ^*(\sqrt{n})$ under the action of $G$.}
In \cite{Riaz1} M. A. Malik and M. Riaz obtained $G$-subsets of $\QQ^*(\sqrt{n})$ by the help of $G^{**}$ acting on $\QQ^*(\sqrt{n})$, when $n\equiv0~(mod~p)$. We also proved for $h=2k+1\geq3$  there are exactly two $G$-orbits of $\QQ^*(\sqrt{2^h})$ namely $ (2^k\sqrt{2})^G$ and $(\frac{2^k\sqrt{2}}{-1})^G$. In \cite{R4} M. A. Malik et al. obtained $G$-subsets of $\QQ^*(\sqrt{n})$ when $n\equiv0~(mod~p^k)$ and  $n\equiv0~(mod~2p^k)$. Thus it becomes interesting to know about more $G$-subsets and to know about the exact number of $G$-orbits under the action of the modular group $G$. In the present studies, we use coset diagrams and modular arithmetic to determine the exact number of $G$-orbits of $\QQ^*(\sqrt{p^k})$, $\QQ^*(\sqrt{2p^k})$,  $\QQ^*(\sqrt{2^2p^k})$, and in general $\QQ^*(\sqrt{2^lp^k})$, for all $l\geq3$ and $k=2h+1\geq3$,   for each odd prime $p$.\\
\\
\textbf{Theorem 2.1}\\
If $p\equiv1~(mod~4)$ and $k=2h+1\geq3$ then there are exactly two $G$-orbits of $\QQ^*(\sqrt{p^k})$ namely
$(\frac{p^h\sqrt{p}}{1})^G$ and $(\frac{1+p^h\sqrt{p}}{2})^G$.\\
\textbf{Proof.}\\
If $p\equiv1~(mod~4)$ then by \cite{R12}, $(\frac{\sqrt{p}}{1})^G$ and $(\frac{1+\sqrt{p}}{2})^G$ are two
distinct $G$-orbits of $\QQ^*(\sqrt{p})$. Then by \cite{Asim} and \cite{R4} the classes $[a,b,c]~(mod p^k$)
with $b$ or $c$ quadratic residues of
$p^k$ lie in the orbit $(\frac{p^h\sqrt{p}}{1})^G$ and similarly  the classes $[a,b,c]~ (mod p^k$),
with $b$ or $c$ quadratic
non-residues of $p^k$ lie in the orbit $(\frac{1+p^h\sqrt{p}}{2})^G$. If $p\equiv1~(mod~4)$
then by \cite{R12},  \cite{R5} $\frac{p^h\sqrt{p}}{-1} \in (\frac{p^h\sqrt{p}}{1})^G$ and none
of $\frac{p^h\sqrt{p}}{1}$ and $\frac{p^h\sqrt{p}}{-1}$ is contained in $(\frac{1+p^h\sqrt{p}}{2})^G$.
 Thus it is clear that $(\frac{p^h\sqrt{p}}{1})^G$ lies in the $A_1$ and $(\frac{1+p^h\sqrt{p}}{2})^G$ lies
  in $A_2$. Moreover the  $G$-subsets $A_1$ and $A_2$ are transitive.\\
In the closed path lying in the orbit $(\frac{p^h\sqrt{p}}{1})^G$,
the transformation $g$ given by
$$(yx)^{m_1}(y^2x)^{n_1}(yx)^{m_2}(y^2x)^{n_2}...(yx)^{m_k}(y^2x)^{n_k}$$  fixes
$p^h\sqrt{p}$, that is $g(k)=k$, and so gives the quadratic equation
$k^2-p^h\sqrt{p}=0$, the zeros, $\pm {p^h\sqrt{p}}$, of this
equation are fixed points of the transformations $g$. Similarly in
the closed path lying in the orbit $(\frac{1+p^h\sqrt{p}}{2})^G$,
the transformation $g'$ given by
$$(yx)^{m'_1}(y^2x)^{n'_1}(yx)^{m'_2}(y^2x)^{n'_2}...(yx)^{m'_k}(y^2x)^{n'_k}$$
  fixes
$\frac{1+p^h\sqrt{p}}{2}$, this proves the result. \quad\quad$\Box$\\
\\
\textbf{Example 2.2}\\
There are exactly two $G$-orbits of $\QQ^*(\sqrt{5^5})$  namely $\pm
5^2(\sqrt{5})$, In the closed path lying in the orbit
$(5^2\sqrt{5})^G$, the transformation
$$(yx)^{22}(y^2x)^{5}(yx)^{1}(y^2x)^{1}(yx)^{5}(yx)^{22}(y^2x)^{5}(yx)^{1}(y^2x)^{1}(yx)^{5}$$  fixes
$5^2\sqrt{5}$, that is $g(k)=k$, and so gives the quadratic equation
$k^2-5^2\sqrt{5}=0$, the zeros, $\pm {5^2\sqrt{5}}$, of this
equation are fixed points of the transformations $g$. Similarly in
the closed path lying in the orbit $(\frac{1+5^2\sqrt{5}}{2})^G$,
the transformation $(yx)^{5}(y^2x)^{11}(yx)^{6}$  fixes
$\frac{1+5^2\sqrt{5}}{2}$.  \quad\quad $\Box$\\
\textbf{Theorem 2.3}\\
If $p\equiv3~(mod~4)$ and $k=2h+1\geq3$ then there are exactly two
$G$-orbits of $\QQ^*(\sqrt{p^k})$ namely
$(\frac{p^h\sqrt{p}}{1})^G$
and $(\frac{p^h\sqrt{p}}{-1})^G$, for each odd prime $p$.\\
\textbf{Proof.}\\
If $p\equiv3~(mod~4)$ then by \cite{R5} $(\frac{\sqrt{p}}{1})^G$
and $(\frac{\sqrt{p}}{-1})^G$ are two distinct $G$-orbits of $\QQ^*(\sqrt{p})$. By \cite{R4} it is clear that
$(\frac{p^h\sqrt{p}}{1})^G$ lies in the $A_1$ and $(\frac{p^h\sqrt{p}}{-1})^G$ lies in $A_2$.
Moreover the  $G$-subsets $A_1$ and $A_2$ are transitive.\\
In the closed path lying in the orbit $(p^h\sqrt{p})^G$, the
transformation $j$ given by
$$(yx)^{s_1}(y^2x)^{t_1}(yx)^{s_2}(y^2x)^{t_2}...(yx)^{s_k}(y^2x)^{t_k}$$
fixes $\frac{p^h\sqrt{p}}{1}$, that is $j(k)=k$, and so gives the
quadratic equation $k^2-p^h\sqrt{p}=0$, the zeros, $\pm
{p^h\sqrt{p}}$, of this equation are fixed points of the
transformations $j$. Similarly in the closed path lying in the orbit
$(\frac{p^h\sqrt{p}}{-1})^G$, the transformation $j'$ given by
$$(yx)^{s'_1}(y^2x)^{t'_1}(yx)^{s'_2}(y^2x)^{t'_2}...(yx)^{s'_k}(y^2x)^{t'_k}$$
fixes $\frac{p^h\sqrt{p}}{-1}$, this proves the result. \quad\quad$\Box$\\
\\
\textbf{Example 2.4}\\
There are exactly two $G$-orbits of $\QQ^*(\sqrt{3^5})$  namely
$(\frac{3^2\sqrt{3}}{ \pm1})^G$, In the closed path lying in the
orbit $(\frac{3^2\sqrt{3}}{1})^G$, the transformation
$$(yx)^{15}(y^2x)^{1}(yx)^{1}(y^2x)^{2}(yx)^{3}(y^2x)^{15}(yx)^{3}(y^2x)^{2}(yx)^{1}(y^2x)^{1}(yx)^{15}$$  fixes
$3^2\sqrt{3}$, that is\\
$((yx)^{15}(y^2x)^{1}(yx)^{1}(y^2x)^{2}(yx)^{3}(y^2x)^{15}(yx)^{3}(y^2x)^{2}(yx)^{1}(y^2x)^{1}(yx)^{15})(k)=k$,\\
and so gives the quadratic equation $k^2-3^2\sqrt{3}=0$, the zeros,
$\pm {3^2\sqrt{3}}$, of this equation are fixed points of the
transformations $g$. Similarly in the closed path lying in the orbit
$(\frac{3^2\sqrt{3}}{-1})^G$, the transformation\\
$(yx)^{15}(y^2x)^{1}(yx)^{1}(y^2x)^{2}(yx)^{3}(y^2x)^{15}(yx)^{3}(y^2x)^{2}(yx)^{1}(y^2x)^{1}(yx)^{15}$\\
fixes
$(\frac{3^2\sqrt{3}}{-1})^G$.  \quad\quad $\Box$\\
\textbf{Theorem 2.5}\\
If $p\equiv1~(mod~4)$ and $k=2h+1\geq3$ then there are exactly two
$G$-orbits of $\QQ^*(\sqrt{2p^k})$ namely
$(\frac{p^h\sqrt{2p}}{1})^G$ and $(\frac{1+p^h\sqrt{2p}}{2})^G$.\\
\textbf{Proof.}\\
Since $1$ is the only quadratic residue of $2$ and there is no
quadratic non-residue of $2$. Thus  the quadratic
residues and quadratic non residues of $p^k$ and $2p^k$ are the
same. Then by \cite{Asim} and \cite{R4} the classes $[a,b,c]$ (modulo $p^k$) with $b$ or $c$
quadratic residues of $p^k$ lie in the orbit $(p^h\sqrt{p})^G$ and
similarly the classes $[a,b,c]$ (modulo $p^k$), with $b$ or $c$
quadratic non-residues of $p^k$ lie in the orbit
$(\frac{1+p^h\sqrt{2p}}{2})^G$. If $p\equiv1~(mod~4)$ then by
\cite{R12} $\frac{p^h\sqrt{2p}}{-1} \in (\frac{p^h\sqrt{2p}}{1})^G$ and none of $\frac{p^h\sqrt{2p}}{1}$
and $\frac{p^h\sqrt{2p}}{-1}$ is contained in
$(\frac{1+p^h\sqrt{p}}{2})^G$. Thus it clear that
$(\frac{p^h\sqrt{2p}}{1})^G$ and $(\frac{1+p^h\sqrt{2p}}{2})^G$
are distinct orbits.\\
By \cite{R4} it is clear that $(\frac{p^h\sqrt{p}}{1})^G$ lies in
the $A_1$ and  $(\frac{1+p^h\sqrt{p}}{2})^G$ lies in $A_2$.\\
In the closed path lying in the orbit $(\frac{p^h\sqrt{2p}}{1})^G$,
the transformation $g$ given by
$$(yx)^{m_1}(y^2x)^{n_1}(yx)^{m_2}(y^2x)^{n_2}...(yx)^{m_k}(y^2x)^{n_k}$$  fixes
$p^h\sqrt{p}$, that is $g(k)=k$, and so gives the quadratic equation
$k^2-p^h\sqrt{2p}=0$, the zeros, $\pm {p^h\sqrt{2p}}$, of this
equation are fixed points of the transformations $g$. Similarly in
the closed path lying in the orbit $(\frac{1+p^h\sqrt{2p}}{2})^G$,
the transformation
 $g'$ given by
$$(yx)^{m'_1}(y^2x)^{n'_1}(yx)^{m'_2}(y^2x)^{n'_2}...(yx)^{m'_k}(y^2x)^{n'_k}$$
  fixes  fixes
$\frac{p^h\sqrt{2p}}{-1}$, this proves the result. \quad\quad$\Box$\\
\textbf{Theorem 2.6}\\
If $p\equiv3~(mod~4)$ and  $k=2h+1\geq3$ then there are exactly
two $G$-orbits of $\QQ^*(\sqrt{2p^k})$ namely
$(\frac{p^h\sqrt{2p}}{1})^G$ and $(\frac{p^h\sqrt{2p}}{-1})^G$, for each odd prime $p$.\\
\textbf{Proof.}\\
  Since $1$ is the only quadratic residue of $2$ and there is no
quadratic non-residue of $2$. Thus the quadratic residues and quadratic non residues of $p^k$ and $2p^k$ are the
same. The classes $[a,b,c]$ (modulo $2p^k$) with $b$ or $c$
quadratic residues of $2p^k$ lie in the orbit $(p^h\sqrt{2p})^G$
and similarly
 the classes $[a,b,c]$ (modulo $2p^k$), with $b$ or $c$ quadratic
non-residues of $2p^k$ lie in the orbit
$(\frac{p^h\sqrt{2p}}{-1})^G$. If $p\equiv3~(mod~4)$ then by
\cite{R5} $\frac{p^h\sqrt{p}}{1}$ and $(\frac{p^h\sqrt{p}}{-1})^G$
are distinct orbits.\\
By \cite{R4} it is clear that $(\frac{p^h\sqrt{p}}{1})^G$ lies in
the $A_1$ and  $(\frac{p^h\sqrt{p}}{-1})^G$ lies in $A_2$.\\
In the closed path lying in the orbit $(p^h\sqrt{2p})^G$, the
transformation $j$ given by
$$(yx)^{s_1}(y^2x)^{t_1}(yx)^{s_2}(y^2x)^{t_2}...(yx)^{s_k}(y^2x)^{t_k}$$
fixes $\frac{p^h\sqrt{p}}{1}$, that is $j(k)=k$, and so gives the
quadratic equation $k^2-p^h\sqrt{2p}=0$, the zeros, $\pm
{p^h\sqrt{2p}}$, of this equation are fixed points of the
transformations $g$. Similarly in the closed path lying in the orbit
$(\frac{p^h\sqrt{2p}}{-1})^G$, the transformation $j'$ given by
$$(yx)^{s'_1}(y^2x)^{t'_1}(yx)^{s'_2}(y^2x)^{t'_2}...(yx)^{s'_k}(y^2x)^{t'_k}$$  fixes
$\frac{p^h\sqrt{2p}}{-1}$, This proves the result. \quad\quad$\Box$\\
\textbf{Theorem 2.7}\\
If $p\equiv1~(mod~4)$ and $k=2h+1\geq3$ then there are exactly two
$G$-orbits of $\QQ^*(\sqrt{2^2p^k})$ namely
$(\frac{2p^h\sqrt{p}}{1})^G$ and $(\frac{1+2p^h\sqrt{p}}{2})^G$.\\
\textbf{Proof.}\\
The proof is analogous to the proof of Theorem 2.5.\\
\textbf{Theorem 2.8}\\
If $p\equiv3~(mod~4)$ and  $k=2h+1\geq3$ then there are exactly
two $G$-orbits of $\QQ^*(\sqrt{2^2p^k})$ namely
$(\frac{p^h\sqrt{2p}}{1})^G$
and $(\frac{p^h\sqrt{2p}}{-1})^G$, for each odd prime $p$.\\
\textbf{Proof.}\\
The proof is analogous to the proof of Theorem 2.6.\\
\\
The results becomes interesting in modulo $8$ and we observe that there are exactly four  $G$-orbits when  $2^lp^k \equiv 0(mod 8)$,
$l\geq3$.\\
\\
\textbf{Theorem 2.9}\\
Let $k=2h+1\geq3$ and  $l\geq3$, Then there are exactly four
$G$-orbits of $\QQ^*(\sqrt{2^lp^k})$ namely
$(\frac{p^h\sqrt{2^lp}}{1})^G$, $(\frac{p^h\sqrt{2^lp}}{-1})^G$,
$(\frac{1+p^h\sqrt{2^lp}}{3})^G$
 and  $(\frac{-1+p^h\sqrt{2^lp}}{-3})^G$ for each odd prime $p$.\\
\textbf{Proof.}\\
We know by Theorem 1.1 \cite{Asim} that if $n\equiv 0(mod~2^3)$ then
$A=\{\alpha\in Q^*(\sqrt{n}):~\textmd{either}~c\equiv
i(mod~2^3)~\textmd{or
 if}~c\equiv0~\textmd{or}~4(mod~2^2)~\textmd{then}~b\equiv
i(mod~2^3)\}$ is a $G$-subset of $Q^*(\sqrt{n})$ for each $i=1,3,5$ and $7$,
we represent these $G$-subsets by $A_1$, $A_2$, $A_3$ and $A_4$.\\
Thus it is easy to see that the orbits
$(\frac{p^h\sqrt{2^lp}}{1})^G$, $(\frac{p^h\sqrt{2^lp}}{-1})^G$,
$(\frac{1+p^h\sqrt{2^lp}}{3})^G$
 and  $(\frac{1+p^h\sqrt{2^lp}}{-3})^G$ are contained in $A_1$,
 $A_2$, $A_3$ and $A_4$ respectively. Moreover it is clear that
 $G$-subsets  given by $A_1$,  $A_2$, $A_3$ and $A_4$ are
 transitive.\\
 In the closed path lying in the orbits $(\frac{p^h\sqrt{2^lp}}{\pm1})^G$, the
transformation
$$(yx)^{u_1}(y^2x)^{v_1}(yx)^{u_2}(y^2x)^{v_2}...(yx)^{u_k}(y^2x)^{v_k}$$ fixes
$\frac{p^h\sqrt{2p}}{\pm1}$. In the closed path lying in the orbits
$(\frac{1+p^h\sqrt{2^lp}}{3})^G$
 and  $(\frac{-1+p^h\sqrt{2^lp}}{-3})^G$, the transformation
$$(yx)^{u'_1}(y^2x)^{v'_1}(yx)^{u'_2}(y^2x)^{v'_2}...(yx)^{u'_k}(y^2x)^{v'_k}$$   fixes
$\frac{-1+p^h\sqrt{2^lp}}{\pm3}$, This proves the result. \quad\quad$\Box$\\
\\
\textbf{Example 2.10}\\
There are exactly four $G$-orbits of $\QQ^*(\sqrt{2^5.3^7})$
namely $(\frac{2^2.3^3\sqrt{2.3}}{1})^G$, $(\frac{2^2.3^3\sqrt{2.3}}{1})^G$, $(\frac{2^2.3^3\sqrt{2.3}}{3})^G$
 and  $(\frac{2^2.3^3\sqrt{2.3}}{-3})^G$.\\
Similarly there are exactly four $G$-orbits of $\QQ^*(\sqrt{2^6.3^7})$ namely $(\frac{2^3.3^3\sqrt{3}}{1})^G$,
$(\frac{2^3.3^3\sqrt{3}}{1})^G$, $(\frac{2^3.3^3\sqrt{3}}{3})^G$  and  $(\frac{2^3.3^3\sqrt{3}}{-3})^G$.\\
\\
\textbf{Conclusion} The action of the modular group $PSL(2,\mathbb {Z})$ on $\QQ^*(\sqrt{p^k})$,
 $k=2h+1\geq3$ is intransitive. $\QQ^*(\sqrt{2^lp^k})$ splits into exactly two $G$-orbits for each
  $l\geq3$ and $k=2h+1 \geq3$, Moreover $\QQ^*(\sqrt{2^lp^k})$ splits into exactly four $G$-orbits
   for each  $l\geq3$ and $k=2h+1 \geq3$ or $ 2^lp^k\equiv 0(mod~ 8)$.\\


\vspace{0.5cm}

\begin{thebibliography}{99}


\bibitem{Anna} A. Torstensson: Coset diagrams in the study of finitely presented groups with an
application to quotients of the modular group, Journal of
commutative algebra 2(4), 501-514, 2010.


\bibitem{Gshabbir} G. Shabbir, T. A. Khan: Proper projective symmetry in some well known conformally flat
space-times,  U. P. B. Scientific Bulletin., Series A 70(1), 25-34,
 2008.


\bibitem{Asim} M. A. Malik, M. A. Zafar: Real Quadratic
Irrational  Numbers and Modular Group Action, Southeast Asian
Bulletin of Mathematics  35(3), 439-445, 2011.

\bibitem{Riaz1} M. A. Malik, M. Riaz: $G$-subsets and $G$-orbits of $\QQ^*(\sqrt{n})$
under action of the Modular Group,  Punjab University Journal of
Mathematics 43, 75-84, 2011.

\bibitem{R4} M. A. Malik, M. Qayyum: Quadratic Irrational Numbers in Modulo $p^k$
 and the Modular Group (Submitted, 2011).

\bibitem{R11} M. A. Malik, S. M. Husnine and A. Majeed:  Modular Group Action on Certain Quadratic
Fields,  Punjab University Journal of Mathematics 28, 47-68, 1995.

 \bibitem{R12} M. A. Malik, S. M. Husnine and A. Majeed:  The orbits of $\QQ^*(\sqrt{p})$, $p=2~\textmd{or}\equiv1~(mod~4)$,
  Under the Action of the Modular Group, Punjab University Journal of Mathematics 33,
  37-50, 2000.

\bibitem{R5} M. A. Malik, S. M. Husnine and A. Majeed:  The orbits of $\QQ^*(\sqrt{p})$, $p\equiv3~(mod~4)$,
Under the Action of the Modular Group   $G=\langle
x,y:x^2=y^3=1\rangle$, Punjab University Journal of Mathematics 36,
1-13, 2003-2004.


\bibitem{R15} M. H. Millington: Subgroups of the classical modular
group,  J. London Math. Soc. 1, 351-357, 1969.

\bibitem{Novotny} P. Novotny, J. Hrivnak: On orbits of the ring $ \mathbb{Z}^m_n$ under action of the group $SL(m,\mathbb
{Z}_n)$, Acta Polytechnica 45(5), 39-43, 2005.

\bibitem{R10} Q. Mushtaq: On word structure of the Modular Group over finite and real quadratic
fields, Discrete Mathematics 179, 145-154, 1998.

\bibitem{R14}S. Anis, Q. Mushtaq: The Number of Subgroups of $PSL(2,\mathbb {Z})$ when acting on
$F_p\cup\{\infty\}$, Communication in Algebra 36, 4276-4283, 2008.



\end{thebibliography}
\end{document}